\documentstyle[amssymb,12pt]{article}

\begin{document}

\author{Mircea Crasmareanu\thanks{%
Partially supported by the CEex Grant 05-D11-84}}
\title{Last multipliers as autonomous solutions of the Liouville equation of
transport}
\date{}
\maketitle

\begin{abstract}
Using the characterization of last multipliers as solutions of the
Liouville's transport equation, new results are given in this
approach of ODE by providing several new characterizations, e.g.
in terms of Witten and Marsden differentials or adjoint vector
field. Applications to Hamiltonian vector fields on Poisson
manifolds and vector fields on Riemannian manifolds are presented.
In Poisson case, the unimodular bracket considerably simplifies
computations while, in the Riemannian framework, a Helmholtz type
decomposition yields remarkable examples: one is the quadratic
porous medium equation, the second (the autonomous version of the
previous) produces harmonic square functions, while the third
refers to the gradient of the distance function with respect to a
two dimensional rotationally symmetric metric. A final example
relates the solutions of Helmholtz (particularly Laplace) equation
to provide a last multiplier for a gradient vector field. A
connection of our subject with gas dynamics in Riemannian setting
is pointed at the end.
\end{abstract}



\vspace{.3cm}

\section*{Introduction}

In January 1838, Joseph Liouville(1809-1882) published a note (\cite{j:l})
on the time-dependence of the Jacobian of the "transformation" exerted by
the solution of an ODE on its initial condition. In modern language, if $%
A=A(x)$ is the vector field corresponding to the given ODE and $m=m(t,x)$ is
a smooth function (depending also of the time $t$), then the main equation
of the cited paper is:
$$
\frac{dm}{dt}+m\cdot divA=0 \eqno(LE)
$$
called, by then, the {\it Liouville equation}. Some authors use the name
{\it generalized Liouville equation} (\cite{g:e}) but we prefer to name it
the {\it Liouville equation of transport} (or {\it of continuity}). This
equation is a main tool in statistical mechanics where a solution is called
a {\it probability density function} (\cite{g:u}).

The notion of {\it last multiplier}, introduced by Carl Gustav Jacob Jacobi
(1804-1851) around 1844, was treated in detailed in ''Vorlesugen \"{u}ber
Dynamik'', edited by R. F. A. Clebsch in Berlin in 1866. So that, sometimes
it has been used under the name of ''Jacobi multiplier''. Since then, this
tool for understanding ODE was intensively studied by mathematicians in the
usual Euclidean space ${I\!\!R}^{n}$, cf. the bibliography of \cite{m:c1},
\cite{n:c1}-\cite{n:c4}. For all those interested in historical aspects an
excellent survey can be found in \cite{b:e}.

The aim of the present paper is to show that last multipliers are
exactly the autonomous, i.e. time-independent, solutions of LE and
to discuss some results of this useful theory extended to
differentiable manifolds. Our study has been inspired by the results
presented in \cite{o:z} using the calculus on manifolds especially
the Lie derivative. Since the Poisson and Riemannian geometries are
the most frequently used frameworks, a Poisson bracket and a
Riemannian metric are added and cases yielding last multipliers are
characterized in terms of unimodular Poisson brackets and
respectively, harmonic functions.

The paper is structured as follows. The first section reviews the
definition of last multipliers and some previous results. New
characterizations in terms of de Rham cohomology and other types
of differentials than the usual exterior derivative, namely Witten
and Marsden, are given. There follows that the last multipliers
are exactly the first integrals of the adjoint vector field. For a
fixed smooth function $m$, the set of vector fields admitting $m$
as last multiplier is a Lie subalgebra of the Lie algebra of
vector fields.

In the next section the Poisson framework is discussed, on showing
that some simplifications are possible for the unimodular case.
For example, in unimodular setting, the set of functions $f$ which
are last multiplier for exactly the Hamiltonian vector field
generated by $f$ is a Poisson subalgebra. Local expressions for
the main results of this section are provided in terms of the
bivector $\pi $ defining the Poisson bracket. Two examples are
given: one with respect to dimension two when as a last multiplier
is obtained exactly the function defining the bivector, and the
second related to Lie-Poisson brackets when the condition to be
last multiplier is expressed in terms of structural constants of
given Lie algebra. Again, the two dimensional case is studied in
detail and for non-vanishing structure constants, an affine
function is obtained as a last multiplier.

The last section, is devoted to the Riemannian manifolds and again,
some characterizations are given in terms of vanishing of some
associated differential operators, e.g. the codifferential. Assuming
a Helmholtz type decomposition, several examples are given: the
first is related to a parabolic equation of porous medium type and
the second yields harmonic square functions. As to the first
example, one should notice the well-known relationship between the
heat equation (in our case, a slight generalization) and the general
method of multipliers; see the examples from \cite[p. 364]{r:j}. The
third example is devoted to the distance function on a two
dimensional rotationally symmetric Riemannian manifold, while the
final example is concerned with solutions of Helmholtz (particularly
Laplace) equation, for obtaining a last multiplier for a gradient
vector field.

The last section deals with an application of our framework to gas
dy\-na\-mics on Riemannian manifolds but in order to not enlarge our
paper too much we invite the reader to see details of physical
nature in Sibner's paper \cite{s:s}. More precisely, a main result
from \cite{s:s} is rephrased in terms of last multipliers.

{\bf Acknowledgments} The author expresses his thanks to Ioan
Bucataru, Marian-Ioan Munteanu, Zbiegnew Oziewicz and Izu Vaisman
for several useful remarks. Also, David Bao and an anonymous referee
suggest some important improvements.

\section{General facts on last multipliers}

Let $M$ be a real, smooth, $n$-dimensional manifold, $C^{\infty }\left(
M\right) $ the algebra of smooth real functions on $M$, ${\cal X}\left(
M\right) $ the Lie algebra of vector fields and $\Lambda ^{k}\left( M\right)
$ the $C^{\infty }\left( M\right) $-module of $k$-differential forms, $0\leq
k\leq n$. Assume that $M$ is orientable with the fixed volume form $V\in
\Lambda ^{n}\left( M\right) $.

Let:
\[
\stackrel{.}{x}^{i}\left( t\right) =A^{i}\left( x^{1}\left( t\right) ,\ldots
,x^{n}\left( t\right) \right) ,1\leq i\leq n
\]
an ODE system on $M$ defined by the vector field $A\in {\cal X}\left(
M\right) ,A=\left( A^{i}\right) _{1\leq i\leq n}$ and let us consider the $%
\left( n-1\right) $-form $\Omega =i_{A}V\in \Lambda ^{n-1}\left( M\right) $.

\bigskip

{\bf Definition 1.1}(\cite[p. 107]{f:f}, \cite[p. 428]{o:z}) The function $%
m\in C^{\infty }\left( M\right) $ is called a {\it last multiplier} of the
ODE system generated by $A$, ({\it last multiplier }of $A$, for short) if $%
m\Omega $ is closed:
$$
d\left( m\Omega \right) :=\left( dm\right) \wedge \Omega +md\Omega =0.\eqno%
\left( 1.1\right)
$$

For example, in dimension $2$, the notions of last multiplier and
integrating factor are identical and Sophus Lie suggests a method to
associate a last multiplier to every symmetry vector field of $A$ (Theorem
1.1 in \cite[p. 752]{h:s}). Lie method is extended to any dimension in \cite
{o:z}.

If the $\left( n-1\right)^{th} $ de Rham cohomology space of $M$ is zero, $%
H^{n-1}\left( M\right) =0$, there follows that $m$ is a last
multiplier iff there exists $\alpha \in \Lambda ^{n-2}\left(
M\right) $ so that $m\Omega =d\alpha $. Other characterizations
can be obtained in terms of Witten's differential \cite{w:e} and
Marsden's differential \cite[p. 220]{m:a}. If $f\in C^{\infty
}\left( M\right) $ and $t\geq 0$, Witten deformation of the usual
differential $d_{tf}:\Lambda ^{\ast }\left( M\right) \rightarrow
\Lambda ^{\ast +1}\left( M\right) $ is defined by:
\[
d_{tf}=e^{-tf}de^{tf}
\]
which means \cite{w:e}:
\[
d_{tf}\left( \omega \right) =tdf\wedge \omega +d\omega .
\]
Hence, $m$ is a last multiplier if and only if:
\[
d_{m}\Omega =\left( 1-m\right) d\Omega
\]
i.e. $\Omega $ belongs to the kernel of the differential operator $%
d_{m}+\left( m-1\right) d:\Lambda ^{n-1}\left( M\right) \rightarrow \Lambda
^{n}\left( M\right) $. Marsden differential is $d^{f}:\Lambda ^{\ast }\left(
M\right) \rightarrow \Lambda ^{\ast +1}\left( M\right) $ defined by:
\[
d^{f}\left( \omega \right) =\frac{1}{f}d\left( f\omega \right)
\]
and $m$ is a last multiplier iff $\Omega $ is $d^{m}$-closed.

The following characterization of last multipliers will be useful:

{\bf Lemma 1.2}(\cite[p. 428]{o:z}) (i) $m\in C^{\infty }\left( M\right) $
{\it is a last multiplier for }$A$ {\it if and only if}:
$$
A\left( m\right) +m\cdot divA=0\eqno\left( 1.2\right)
$$
{\it where $divA$ is the divergence of $A$ with respect to volume form $V$.}

(ii) {\it Let $0\neq h\in C^{\infty }\left( M\right) $ such that:
$$
L_{A}h:=A\left( h\right) =\left( divA\right) \cdot h\eqno\left( 1.3\right)
$$
Then $m=h^{-1}$ is a last multiplier for $A$}.

{\bf Remarks} {\bf 1.3 } (i) Equation $(1.2)$ is exactly the
time-independent version of LE from the Introduction. So that, the promised
relationship between LE and last multipliers is obtained. An important
feature of equation $\left( 1.2\right) $ is that it does not always admit
solutions cf. \cite[p. 269]{g:fc}.

(ii) In the terminology of \cite[p. 89]{b:e}, a function {\it h} satisfying
(1.3) is called an {\it inverse multiplier}.

(iii) A first result given by $\left( 1.2\right) $ is the
characterization of last multipliers for divergence-free vector
fields: $m\in C^{\infty }\left( M\right) $ {\it is a last multiplier
for the divergenceless vector field }$A$ {\it iff }$m$ {\it is a
first integral of} $A$. The importance of this result is shown by
the fact that three remarkable classes of divergence-free vector
fields are provided by: Killing vector fields in Riemannian
geometry, Hamiltonian vector fields in symplectic geometry and Reeb
vector fields in contact geometry. Also, there are many equations of
mathematical physics corresponding to the vector fields without
divergence.

(iv) For the general case, namely $A$ is not divergenceless, there
is a strong connection between first integrals and last
multipliers as well. Namely, from properties of Lie derivative,
the ratio of two last multipliers is a first integral and
conversely, the product between a first integral and a last
multiplier is a last multiplier. So, denoting $FInt(A)$ the set of
first integrals of $A$, since $FInt(A)$ is a subalgebra in
$C^{\infty }(M)$ it results that the set of last multipliers for
$A$ is a $FInt(A)$-module.

(v) Recalling formula:
$$
div\left( fX\right) =X\left( f\right) +fdivX \eqno(1.4)
$$
there follows that $m$ is a last multiplier for $A$ if and only if
the vector field $mA$ is without divergence i.e. $div\left(
mA\right) =0$. Thus, the set of last multipliers is a "measure of
how far away" is $A$ from being divergence-free.

(vi) To the vector field $A$ one may associate an {\it adjoint
vector field} $A^{\ast }$, acting on functions in the following
manner, \cite[p. 129]{m:t}:
\[
A^{\ast }\left( m\right) =-A\left( m\right) -mdivA.
\]
Then, another simple characterization is: $m$ is a last multiplier
for $A$ iff $m$ is a first integral of the adjoint $A^{\ast }$. An
important consequence results: the set of last multipliers is a
subalgebra in $C^{\infty }\left( M\right) $.

An important structure generated by a last multiplier is given by:

{\bf Proposition 1.4} {\it Let }$m\in C^{\infty }\left( M\right) ${\it \ be
fixed. The set of vector fields admitting }$m${\it \ as last multiplier is a
Lie subalgebra in} ${\cal X}\left( M\right) $.

{\bf Proof} Let $X$ and $Y$ be vector fields with the required property.
Since \cite[p. 123]{m:r}:
\[
div\left[ X,Y\right] =X\left( divY\right) -Y\left( divX\right)
\]
one has:
\[
\left[ X,Y\right] \left( m\right) +mdiv\left[ X,Y\right] =\left( X\left(
Y\left( m\right) \right) +mX\left( divY\right) \right) -\left( Y\left(
X\left( m\right) \right) +mY\left( divX\right) \right) =
\]
\[
=\left( -divY\cdot X\left( m\right) \right) -\left( -divX\cdot Y\left(
m\right) \right) =divY\cdot mdivX-divX\cdot mdivY=0.\qquad \square
\]

\section{Last multipliers on Poisson manifolds}

Let us assume that $M$ is endowed with a Poisson bracket $\{,\}$. Let $f\in
C^{\infty }\left( M\right) $ and $A_{f}\in {\cal X}\left( M\right) $ be the
associated {\it Hamiltonian vector field} of the {\it Hamiltonian} $f$ cf.
\cite{m:r}. Recall that given the volume form $V$ there exists a unique
vector field $X_{V},$ called the {\it modular vector field}, so that (\cite
{k:l}, \cite{a:w}):
\[
div_{V}A_{f}=X_{V}\left( f\right) .
\]

The triple $\left( M,\{,\},V\right) $ is called (\cite{a:w}) {\it unimodular}
if $X_{V}$ is a Hamiltonian vector field, $A_{\rho}$ of $\rho \in C^{\infty
}\left( M\right) $.

From $\left( 1.2\right) $ there results:
\[
0=A_{f}\left( m\right) +mX_{V}\left( f\right) =-A_{m}\left( f\right)
+mX_{V}\left( f\right)
\]
which means:

{\bf Proposition 2.1} $m$ {\it is a last multiplier of} $A_{f}$ {\it if and
only if }$f${\it \ is a first integral for the vector field} $mX_{V}-A_{m}$.
{\it In the unimodular case, }$m${\it \ is a last multiplier for }$A_{f}$%
{\it \ if and only if }$m\{\rho ,f\}=\{m,f\}${\it .}

Since $f$ is a first integral of $A_{f}$ we get:

{\bf Corollary 2.2} $f$ {\it is a last multiplier for }$A_{f}${\it \ if and
only if }$f${\it \ is a first integral of the vector field} $X_{V}$. {\it In
the unimodular case, }$f${\it \ is a last multiplier for }$A_{f}${\it \ if
and only if }$\{\rho ,f\}=0${\it .}

With the Jacobi and Leibniz formulas the following consequence of Corollary
2.2 may be established:

{\bf Corollary 2.3} {\it Let $\left( M,\{,\}\right) $ be a unimodular Poison
manifold and let }$F${\it \ be the set of smooth functions }$f${\it \ that
are last multipliers of }$A_{f}${\it . Then }$F${\it \ is a Poisson
subalgebra in }$\left( C^{\infty }\left( M\right) ,\cdot ,\{,\}\right) .$

Another important consequence of Proposition 2.1 is:

{\bf Corollary 2.4} {\it If $m$ is a last multiplier of $A_{f}$ and $A_{g}$
then $m$ is a last multiplier of $A_{fg}$. Then, if $m$ is a last multiplier
of $A_{f}$ then $m$ is a last multiplier of $A_{f^{r}}$ for every natural
number $r\geq 1$.}

Let $\left( x^{1},\ldots ,x^{n}\right) $ be a local chart on $M$ such that $%
V=dx^{1}\wedge \ldots \wedge dx^{n}$ and the bivector $\pi $ of $\left(
M,\{,\}\right) $ is: $\pi =\sum\limits_{i<j}\pi ^{ij}\frac{\partial }{%
\partial x^{i}}\wedge \frac{\partial }{\partial x^{j}}$. Denoting $\pi
^{i}=\sum\limits_{j=1}^{n}\frac{\partial \pi ^{ij}}{\partial x^{j}}$ we have
(\cite[Proposition 1, p. 4]{d:f}):
$$
X_{V}=\sum\limits_{i=1}^{n}\pi ^{i}\frac{\partial }{\partial x^{i}}\eqno%
\left( 2.1\right)
$$
and then, Proposition 2.1 and Corollary 2.2 become:

{\bf Proposition 2.5} (i) $m\in C^{\infty }\left( M\right) $ {\it is a last
multiplier for }$A_{f}$ {\it if and only if}:
$$
m\sum\limits_{i=1}^{n}\pi ^{i}\frac{\partial f}{\partial x^{i}}%
=\{m,f\}=\sum\limits_{i<j}\pi ^{ij}\frac{\partial m}{\partial x^{i}}\frac{%
\partial f}{\partial x^{j}}.\eqno(2.2)
$$
(ii) $f\in C^{\infty }\left( M\right) $ {\it is a last multiplier for }$%
A_{f} ${\it \ if and only if}:
$$
\sum\limits_{i=1}^{n}\pi ^{i}\frac{\partial f}{\partial x^{i}}=0.\eqno\left(
2.3\right)
$$

{\bf Examples:}

{\bf 2.1}

After \cite[p. 31]{va:i} the bivector $\pi =h\left( x,y\right)
\frac{\partial }{\partial x}\wedge \frac{\partial }{\partial y}$
defines a
Poisson structure on ${I\!\!R}^{2}$. So, $\pi ^{1}=\frac{\partial h}{\partial y}%
,\pi ^{2}=-\frac{\partial h}{\partial x}$ and then $\left( 2.3\right) $
becomes:
\[
\frac{\partial h}{\partial y}\frac{\partial f}{\partial x}-\frac{\partial h}{%
\partial x}\frac{\partial f}{\partial y}=0
\]
with the obvious solution $f=h$. Therefore, on the Poisson
manifold $\left( {I\!\!R}^{2},\pi \right) $ above the function $h$
is a last multiplier for exactly the Hamiltonian vector field
$A_{h}$.

{\bf 2.2 Lie-Poisson structures}

Let ${\cal G}$ be an $n$-dimensional Lie algebra with a fixed basis $%
B=\{e_{i}\}_{1\leq i\leq n}$ and let $B^{\ast }=\{e^{i}\}$ be the dual basis
on the dual ${\cal G}^{\ast }$. Recall the definition of {\it structure
constants} of ${\cal G}$:
\[
\left[ e_{i},e_{j}\right] =c_{ij}^{k}e_{k}.
\]
Then, on ${\cal G}^{\ast }$ we have the so-called {\it Lie-Poisson
structure} given by (\cite[p. 31]{va:i}):
$$
 \pi ^{ij}\left(
x_{u}e^{u}\right) =c_{ij}^{k}x_{k}.\eqno\left( 2.4\right)
$$
Hence:
$$
\pi ^{i}=\sum\limits_{j=1}^{n}c_{ij}^{j}\eqno\left( 2.5\right)
$$
yields:

{\bf Proposition 2.6} {\it For a Lie-Poisson structure on }${\cal
G}^{\ast }${\it \ provided by the structural constants }$\left(
c_{ij}^{k}\right) ${\it a function} $f\in C^{\infty }\left( {\cal
G}^{\ast }\right) $ {\it is a last multiplier for }$A_{f}${\it \ if
and only if}:
$$
\sum\limits_{i,j=1}^{n}c_{ij}^{j}\frac{\partial f}{\partial x_{i}}=0.\eqno%
\left( 2.6\right)
$$

{\bf Particular case: n=2}

From $\left[ e_{1},e_{1}\right] =[e_{2},e_{2}]=0,\left[ e_{1},e_{2}\right]
=c_{12}^{1}e_{1}+c_{12}^{2}e_{2}$ the last equation reads:
$$
c_{12}^{2}\frac{\partial f}{\partial x_{1}}-c_{12}^{1}\frac{\partial f}{%
\partial x_{2}}=0.\eqno\left( 2.7\right)
$$
For example, if $c_{12}^{1}\cdot c_{12}^{2}\neq 0$, a solution of $\left(
2.7\right) $ is:
$$
f=f\left( x_{1},x_{2}\right) =A(\frac{x_{1}}{c_{12}^{1}}+\frac{x_{2}}{%
c_{12}^{2}})+B\eqno\left( 2.8\right)
$$
with $A,B$ real constants.

\section{The Riemannian case}

Let us suppose that a Riemannian metric $g=<,>$ on $M$ is given; then, there
exists an induced volume form $V_{g}$. Let $\omega \in \Lambda ^{1}\left(
M\right) $ be the $g$-dual of $A$ and $\delta $ the co-derivative operator $%
\delta :\Lambda ^{\ast }\left( M\right) \rightarrow \Lambda ^{\ast -1}\left(
M\right) $. Then:
\[
\left\{
\begin{array}{c}
div_{V_{g}}A=-\delta \omega \\
A\left( m\right) =g^{-1}\left( dm,\omega \right)
\end{array}
\right.
\]
and condition $\left( 1.2\right) $ means:
\[
g^{-1}\left( dm,\omega \right) =m\delta \omega .
\]
Supposing that $m>0$ it follows that $m$ is a last multiplier if and only if
$\omega $ belongs to the kernel of the differential operator: $g^{-1}\left(
d\ln m,\cdot \right) -\delta :\Lambda ^{1}\left( M\right) \rightarrow
\Lambda ^{0}=C^{\infty }\left( M\right) $.

Now, assume that the vector field $A$ admits a Helmholtz type decomposition:
$$
A=X+\nabla u\eqno\left( 3.1\right)
$$
where $X$ is a divergence-free vector field and $u\in C^{\infty }\left( M\right) $%
; for example, if $M$ is compact such decompositions always exist. From $%
div_{V_{g}}\nabla u=\Delta u$, the Laplacian of $u$, and $\nabla u\left(
m\right) =<\nabla u,\nabla m>$ there follows that $\left( 1.2\right) $
becomes:
$$
X\left( u\right) +<\nabla u,\nabla m>+m\cdot \Delta u=0\eqno\left(
3.2\right)
$$

{\bf Example 3.1}

$u$ is a last multiplier of $A=X+\nabla u$ if and only if:
\[
X\left( u\right)=-u\cdot \Delta u-<\nabla u,\nabla u>.
\]
Suppose that $M$ is a cylinder $M=I\times N$ with $I\subseteq {I\!\!R}$ and $%
N$ a $\left( n-1\right) $-manifold; then, for $X=-\frac{1}{2}\frac{\partial
}{\partial t}\in {\cal X}\left( I\right) $ which is divergence-free with
respect to $V=dt\wedge V_N$ with $V_N$ a volume form on $N$, the previous
relation yields:
\[
u_{t}=2\left( u\cdot \Delta u+<\nabla u,\nabla u>\right) .
\]

By the well-known formula (\cite[p. 55]{p:p}):
$$
<\nabla f,\nabla g>=\frac{1}{2}\left( \Delta \left( fg\right) -f\cdot \Delta
g-g\cdot \Delta f\right) \eqno\left( 3.3\right)
$$
the previous equation becomes:
$$
u_{t}=\Delta \left( u^{2}\right) \eqno\left( 3.4\right)
$$
which is a nonlinear parabolic equation of the type of porous medium
equation (\cite{a:r}).

{\bf Example 3.2}

Returning to $\left( 3.1\right) $, suppose that $X=0$. The condition $\left(
3.2\right) $ reads:
$$
m\cdot \Delta u+<\nabla u,\nabla m>=0\eqno\left( 3.5 \right)
$$
which is equivalent, via $(3.3)$ to:
$$
\Delta \left( um\right) +m\cdot \Delta u=u\cdot \Delta m.\eqno\left(
3.6\right)
$$

Adding to $\left( 3.6\right) $ a similar relation with $u$ replaced by $m$
leads to the following conclusion:

{\bf Proposition 3.3} {\it Let }$u,m\in C^{\infty }\left( M\right) ${\it \
such that }$u${\it \ is a last multiplier of }$\nabla m${\it \ and }$m${\it %
\ is a last multiplier of }$\nabla u${\it . Then }$u\cdot m${\it \ is a
harmonic function on }$M$. $u\in C^{\infty }\left( M\right) $ {\it is a last
multiplier of }$A=\nabla u${\it \ if and only if }$u^{2}${\it \ is a
harmonic function on} $M$.

If $M$ is an orientable compact manifold then, from Proposition 3.3, it
follows that $u^{2}$ is a constant which implies that $u$ is a constant,
too. But then $A=\nabla u=0$. Therefore, on an orientable compact manifold,
a function cannot be a last multiplier of its gradient vector field.

If $\left( M,g\right) =({I\!\!R}^n,can)$, there are two classes (with
respect to the sign $\pm $) of radial functions with harmonic square:
\newline
i) $n=2$
\[
u_{\pm }\left( r\right) =\pm \sqrt{C_{1}lnr+C_{2}}
\]
ii) $n=1$ and $n\geq 3$
\[
u_{\pm }(r)=\pm \sqrt{C_1r^{2-n}+C_2}
\]
with $C_{1},C_{2}$ real constants and $r=\sqrt{(x^1)^2+\ldots +(x^n)^2}$.

{\bf Example 3.3. The gradient of the distance function with
respect to a two dimensional rotationally symmetric metric}

Let $M$ be a two dimensional manifold with local coordinates
$(t,\theta )$ endowed with a {\it rotationally symmetric} metric
$g=dt^{2}+\varphi ^{2}(t)d\theta ^{2}$ cf. \cite[p. 11]{p:p}. Let
$u\in C^{\infty }\left( M\right) $, $u\left( t,\theta \right) =t$,
which appears as a distance function with respect to the given
metric. Then, $\nabla u=\frac{\partial }{\partial t}$ and $\Delta
u=\frac{\varphi ^{\prime }\left( t\right) }{\varphi \left(
t\right) }$; equation $\left( 3.5\right) $ is:
\[
m\cdot \frac{\varphi ^{\prime }\left( t\right) }{\varphi \left( t\right) }+%
\frac{\partial m}{\partial t}=0
\]
with solution: $m=m(t)=\varphi ^{-1}\left( t\right) $.

This latter function has a geometric significance: let $T=T\left( t\right) $
be an integral of $m$ i.e. $\frac{dT}{dt}=m=\frac{1}{\varphi \left( t\right)
}$. Then, in the new coordinates $\left( T,\theta \right) $ the given metric
is conformally Euclidean: $g=\varphi ^{2}\left( t\right) \left(
dT^{2}+d\theta ^{2}\right) $ where $t=t\left( T\right) $.

{\bf Example 3.4}

Another way to treat the case $X=0$ is via equation $(1.4)$:
$$
div(m\nabla u)=0 \eqno(3.7)
$$
in which the transformation $v=u\sqrt{m}$ (recall that we search for $m>0$)
is considered. From $div(m\nabla \frac{v}{\sqrt{m}})=0$, it results $div(%
\sqrt{m}\nabla v-v\nabla \sqrt{m})=0$ i.e.:
$$
\sqrt{m}\Delta v=v\Delta \sqrt{m} \eqno(3.8)
$$
which yields:

{\bf Proposition 3.4} {\it Let} $a>0, b$ {\it be solutions of Helmholtz} (%
{\it particularly Laplace}) {\it on the Riemannian manifold} $(M,g)$. {\it %
Then}, $m=a^2$ {\it is a last multiplier for the gradient vector field of
function} $u=\frac{b}{a}$.

\section{Applications to gas dynamics}

Consider again the Riemannian manifold $(M,g)$. Set $m\in C^{\infty }(M)$
and recall according to \cite[p. 62]{s:s}:

{\bf Definition 4.1} A form $\omega $ is said to be: \newline
(i) $m$-{\it coclosed} if $\delta (m\omega )=0$, \newline
(ii) $m$-{\it harmonic} if it is closed and $m$-coclosed.

An important result from the cited paper is:

{\bf Proposition 4.2} {\it If the 1-form} $\omega $ {\it is} $m$-{\it %
harmonic, then, locally,} $\omega =d\phi $, {\it where} $\phi \in C^{\infty
}(M)$ {\it satisfies}:
$$
\delta (md\phi )=0. \eqno(4.1)
$$
In local coordinates $(x^1,...,x^n)$ on $M$ this equation has the
form:
\[
\frac{1}{\sqrt{ det g}}\frac{\partial }{\partial x^i}({\sqrt {det
g}}g^{ij}m\frac{\partial \phi}{\partial x^j})=0
\]
and for a flat $M$ the last equation is exactly the classical gas
dynamics equation conform \cite[p. 63]{s:s}.

But $(4.1)$ is exactly $div_{V_g}(m\nabla \phi)=0$, which means that $m$ is
a last multiplier for the gradient vector field $\nabla \phi $. So, the last
result can be rephrased:

{\bf Proposition 4.3} {\it If the 1-form} $\omega $ {\it is} $m$-{\it %
harmonic, then, locally,} $\omega =d\phi $, {\it where} $\phi \in C^{\infty
}(M)$ {\it with $\nabla \phi$ having $m$ as a last multiplier}.

Due to the local character of previous result for this setting, the
Proposition 1.4 can be improved:

{\bf Proposition 4.4} {\it Let the 1-form} $\omega $ {\it be} $m$-{\it %
harmonic and }$a,b\in C^{\infty }(M)$, {\it such that} $\omega =da=db$.{\it %
\ Then $m$ is a last multiplier for $\nabla a$ and $\nabla b$ but first
integral for the vector field }$\left[ \nabla a,\nabla b\right] $.

{\bf Proof} We have:
\[
\left[ \nabla a,\nabla b\right] (m)=\nabla a\left( m\delta \omega \right)
-\nabla b\left( m\delta \omega \right) =m\left[ \nabla a\left( \delta \omega
\right) -\nabla b\left( \delta \omega \right) \right] .
\]
But $\nabla a\left( \delta \omega \right) =\nabla b\left( \delta \omega
\right) =g^{-1}\left( d\delta \omega ,\omega \right) $. \qquad $\square $

\vspace{.2cm}

\noindent Faculty of Mathematics \newline
University "Al. I. Cuza"\newline
Ia\c si, 700506\newline
Rom\^ania\newline
e-mail: mcrasm@uaic.ro \newline
\newline
\noindent http://www.math.uaic.ro/$\sim$mcrasm


\vspace{0.2cm}

\begin{thebibliography}{99}
\bibitem{a:r}  Aronson, D. G., {\it The porous medium equation}, Nonlinear
diffusion problems (Montecatini Terme, 1985), 1-46, Lecture Notes in Math.,
1224, Springer, Berlin, 1986. MR {\bf 88a}:35130

\bibitem{b:e}  Berrone, L. R.; Giacomini, H., {\it Inverse Jacobi multipliers%
}, Rend. Circ. Mat. Palermo (2) 52(2003), no. 1, 77-130. MR {\bf 2004b}:34067

\bibitem{m:c1}  Crasmareanu, M., {\it Last multipliers theory on manifolds},
Tensor, 66(2005), no. 1, 18-25. MR {\bf 2006e}:34012

\bibitem{g:e}  Ezra, Gregory S., {\it On the statistical mechanics of
non-Hamiltonian systems: the generalized Liouville equation, entropy, and
time-dependent metrics}, J. Math. Chem., 35(2004), No. 1, 29-53. MR {\bf %
2004m}:82073

\bibitem{d:f}  Damianou, Pantelis A.; Fernandes, Rui Loja, {\it Integrable
hierarchies and the modular class}, arXiv: math.DG/0607784, 30 Jul 2006.

\bibitem{f:f}  Flanders, H., {\it Differential forms with applications to
the physical sciences}, Academic Press, 1963. MR {\bf 28} \#5397

\bibitem{g:fc}  Gasc\'{o}n, F. C., {\it Divergence-free vectorfields and
integration via quadrature}, Physics Letters A, 225(1996), 269-273. MR {\bf %
98f}:58015

\bibitem{h:s}  Helgason, S., {\it Invariant differential equations on
homogeneous manifolds}, Bull. AMS, 83(1977), no. 5, 751-774. MR {\bf 56}
\#3579

\bibitem{k:l}  Koszul, J.-L., {\it Crochet de Schouten-Nijenhuis et
cohomologie}, in "\'Elie Cartan et les math\'{e}matiques d'aujourd'hui, The
mathematical heritage of Elie Cartan", Semin. Lyon 1984, Ast\'{e}risque, No.
Hors S\'{e}r. 1985, 257-271 (1985). MR {\bf 88m}:17013

\bibitem{j:l}  Liouville, J., {\it Sur la Th\'eorie de la Variation des
constantes arbitraires}, J. Math. Pures Appl., 3(1838), 342-349.

\bibitem{m:a}  Marsden, Jerrold E., {\it Well-posedness of the equations of
a nonhomogenous perfect fluid}, Comm. Partial Differential Equations,
1(1976), no. 3, 215-230. MR {\bf 53} \#9286.

\bibitem{m:r}  Marsden, Jerrold E.; Ratiu, Tudor S., {\it Introduction to
mechanics and symmetry}, Texts in Applied Math. no. 17, Springer-Verlag,
1994. MR {\bf 2000i}:70002

\bibitem{n:c1}  Nucci, M. C.; Leach, P. G. L., {\it Jacobi's last multiplier
and the complete symmetry group of the Euler-Poinsot system}, J. Nonlinear
Math. Phys. 9(2002), suppl. 2, 110-121. MR {\bf 2003h}:34075

\bibitem{n:c2}  Nucci, M. C.; Leach, P. G. L., {\it Jacobi's last multiplier
and symmetries for the Kepler problem plus a lineal story}, J. Phys. A
37(2004), no. 31, 7743-7753. MR {\bf 2005g}:70012

\bibitem{n:c3}  Nucci, M. C. {\it Jacobi last multiplier and Lie symmetries:
a novel application of an old relationship}, J. Nonlinear Math. Phys.
12(2005), no. 2, 284-304. MR {\bf 2005m}:34017

\bibitem{n:c4}  Nucci, M. C.; Leach, P. G. L., {\it Jacobi's last multiplier
and the complete symmetry group of the Ermakov-Pinney equation}, J.
Nonlinear Math. Phys. 12(2005), no. 2, 305-320. MR {\bf 2006a}:34106

\bibitem{o:z}  Oziewicz, Zbiegniew; Zeni, Jos\'{e} Ricardo R., {\it Ordinary
differential equations: symmetries and last multiplier}, in R. Alamowicz and
B. Fauser(Eds.) ''Clifford Algebras and their Applications in Mathematical
Physics'', vol. 1 (Algebra and Physics), Progr. Phys., vol. 18, Birkhauser,
2000, 425-433. MR {\bf 2001e}:34010

\bibitem{p:p}  Petersen, P., {\it Riemannian geometry}, GTM no. 171,
Springer-Verlag, 1998. MR {\bf 98m}:53001

\bibitem{r:j}  Rauch, J., Review of ''Partial differential equations'' by L.
C. Evans, in Bull. AMS, 37(2000), no. 3, 363-367.

\bibitem{s:s}  Sibner, L. M.; Sibner, R. J., {\it A non-linear Hodge-de-Rham
theorem} Acta Math., 125(1970), 57-73. MR {\bf 43} \#6950

\bibitem{m:t}  Taylor, Michael E., {\it Partial differential equations.
Basic theory}, Texts in Applied Mathematics, no. 23, Springer, N. Y., 1996.
MR {\bf 98b}:35002a

\bibitem{g:u}  \"Unal, Gazanfer, {\it Probability density functions, the
rate of entropy change and symmetries of dynamical systems}, Phys. Lett. A,
233(1997), no. 3, 193-202. MR {\bf 98k}:82107

\bibitem{va:i}  Vaisman, I., {\it Letures on the Geometry of Poisson
Manifolds}, Birkh\"auser, 1994. MR {\bf 95h}:58057

\bibitem{a:w}  Weinstein, Alan, {\it The modular automorphism group of a
Poisson manifold}, J. Geom. Phys. 23(1997), No.3-4, 379-394. MR {\bf 98k}%
:58095

\bibitem{w:e}  Witten, E., {\it Supersymmetry and Morse theory}, J. Diff.
Geom., 17(1982), 661-692. MR {\bf 84b}:58111
\end{thebibliography}
\end{document}